\documentclass[12pt]{article}

\usepackage{amsmath,amsthm,amssymb}
\usepackage{calc,epsfig}
\usepackage{color, graphicx}
\usepackage{latexsym}
\usepackage{amsfonts}
\usepackage{subfigure}

\newcommand{\sig}{{\hat{\Sigma}}_x}
\newcommand{\Xc}{({\bar{X}}_y-{\bar{X}})}
\newcommand{\Xcc}{{\bar{X}}_y-{\bar{X}}}
\newcommand{\Xy}{({\bar{X}}_{y_0}-{\bar{X}})}
\newcommand{\ctsig}{(C_y^T \otimes {\hat{\Sigma}}_x)}
\newcommand{\btsig}{({\hat{C}}_y^T \otimes {\hat{\Sigma}}_x)}
\newcommand{\vecA}{{\rm{vec}}(A)}
\newcommand{\grad}{\nabla}
\newcommand{\R}{\mathbb{R}}
\newcommand{\E}{\mathbb{E}}

\newtheorem{coro}{Corollary}
\newtheorem{prop}{Proposition}

\DeclareMathOperator*{\argmin}{arg\,min}
\newcommand{\CQFD}
{%
\mbox{}%
\nolinebreak%
\hfill%
\rule{2mm}{2mm}%
\medbreak%
\par%
}

\begin{document}

{\title{A Note on Sliced Inverse Regression with Regularizations}}

\author{Caroline Bernard-Michel, Laurent Gardes, St\'ephane Girard$^{\star}$}
\date{\small Laboratoire Jean-Kuntzmann \& INRIA Rh\^one-Alpes, team Mistis,\\
 Inovall\'ee, 655, av. de l'Europe, Montbonnot, 38334 Saint-Ismier cedex, France, \\
($^{\star}$ corresponding author, {\tt Stephane.Girard@inrialpes.fr})
}
\maketitle
\begin{abstract}
\noindent
 In "Li, L. and Yin, X. (2008). Sliced Inverse Regression with Regularizations. {\it{Biometrics}}, 64(1):124--131" 
a ridge SIR estimator is introduced as the solution of a minimization problem
and computed thanks to an alternating least-squares algorithm.
This methodology reveals good performance in practice.
In this note, we focus on the theoretical properties of the estimator.
Is it shown that the minimization problem is degenerated
in the sense that only two situations can occur:
Either the ridge SIR estimator does not exist or it is zero.
 \\


\end{abstract}

\section{Introduction}

Many methods have been developed for inferring the conditional distribution
of an univariate response $Y$ given a predictor $X$ in $\R^p$.
When $p$ is large, sufficient dimension reduction aims at replacing
the predictor $X$ by its projection onto a subspace of smaller dimension
without loss of information on the conditional distribution of
$Y$ given $X$. In this context, the central subspace, denoted by
${\cal S}_{Y|X}$ plays an important role. It is defined as the smallest
subspace such that, conditionally on the projection
of $X$ on ${\cal S}_{Y|X}$, $Y$ and $X$ are independent. In other words, the projection
of $X$ on ${\cal S}_{Y|X}$ contains all the information on $Y$ that is available in 
the predictor $X$. Introducing $d=\dim({\cal S}_{Y|X})$ and 
$A\in\R^{p \times d}$ such that ${\cal{S}}_{Y|X}={\rm{Span}}(A)$,
this property can be rewritten in terms of conditional distribution functions as
\[ F(Y|X) = F(Y|A^TX). \]
The estimation of $A$ has received considerable
attention, and among the proposed methods, 
Sliced Inverse Regression (SIR)~\cite{Li}
seems to be the most popular one.
Let us recall its definition from the minimum discrepancy point of view~\cite{coo04,CookNi}.
Starting from a $n$- sample,
 and denoting by ${\bar{X}}$ the average of $X$, $\sig$ the sample covariance matrix of $X$ and assuming that the response variable $Y$ is partitioned into $h$ non-overlapping slices, the SIR estimator of $A$ is obtained by minimizing
\begin{equation}
\label{SIR}
G(A,C) = \sum_{y=1}^h f_y \left(\Xc - \sig A C_y \right)^T{\sig^{-1}}
 \left(\Xc - \sig A C_y \right)
\end{equation}
where $f_y=n_y/n$, $n_y$ is the number of observations in the $y$th slice, ${\bar{X}}_y$ is the average of $X$ in the $y$th slice and $C=(C_1,\ldots,C_h) \in \R^{d \times h}$. Defining 
$$
\hat\Gamma=\sum_{y=1}^h f_y \Xc \Xc^T,
$$
an estimator of ${\mbox{cov}}(\E(X|Y))$,
the SIR estimator is obtained by computing the eigenvectors of
$\sig^{-1} \hat\Gamma$ associated to the $d$ largest eigenvalues.
It thus
requires the inversion of $\sig$ which is not possible
as soon as $p>n$ or when the predictors are highly correlated.
In order to overcome this problem, it has been proposed
to use the ridge SIR estimator 
(\cite{liyin07}, Definition~1) defined as follows. Let $\tau \geq0$ and
\begin{equation}
\label{SIRreg}
G_{\tau}(A,C) = \sum_{y=1}^h f_y \| \Xc - \sig A C_y \|^2 + \tau \| \vecA \|^2, 
\end{equation}
where ${\rm{vec}}(.)$ is a matrix operator that stacks all columns of the matrix to a single vector. The ridge SIR estimator of the central subspace ${\cal{S}}_{Y|X}$ is ${\rm{Span}}({\hat{A}})$ where
\begin{equation}
\label{argmin}
({\hat{A}},{\hat{C}})=\argmin_{A,C} G_{\tau}(A,C).
\end{equation}
From the practical point of view, an alternating least-squares algorithm
is proposed to solve this optimization problem~\cite{liyin07}.
It revealed good 
performances on simulated and real data. Here, we focus on
the theoretical aspects.
To this end, let us highlight that definition~(\ref{argmin})
assumes the existence of a unique
minimum of $G_\tau$.  In Section~\ref{secexist},
we prove that this is not the case.
In fact, either $\argmin G_\tau=\emptyset$, and thus the ridge SIR estimator
does not exist, or $\argmin G_\tau \subset \{0\}\times  \R^{d \times h}$
and consequently the ridge SIR estimator is zero.
A modification of the criterion~(\ref{SIRreg}) is 
proposed in Section~\ref{secnew} leading 
to the estimator of $A$ proposed in~\cite{Zhong}.
Proofs are postponed to the Appendix.

\section{On the existence of the ridge SIR estimator}
\label{secexist}

Before stating our main result on the existence of the ridge SIR estimator,
remark that $G_\tau$ ($\tau>0$) does not penalize the same way two
proportional matrices $A$ and $\lambda A$, $\lambda\in\R\setminus\{0\}$,
although defining the same central subspace since 
${\rm{Span}}(A)= {\rm{Span}}(\lambda A)$. This lack of invariance
may explain why the ridge SIR estimator is ill-defined as illustrated below.

\begin{prop}
\label{mainprop}
Let $\tau>0$. 
If $\argmin G_\tau\neq \emptyset$ then 
${\hat{A}}$ defined by~(\ref{argmin}) is 
the zero $p\times d$ matrix.  Moreover, 
\begin{equation}
\label{eqprop}
G_{\tau}({\hat{A}},C) = G_{\tau}(0,C) = \sum_{y=1}^h f_y \|\Xcc\|^2,
\end{equation}
for all $C \in \R^{d \times h}$,
\end{prop}
\noindent Since (\ref{eqprop}) does not depend on $C$, it follows that
either $\argmin G_\tau= \emptyset$ or $\argmin G_\tau\subset \{0\}
\times  \R^{d \times h}$.
The following proposition permits to distinguish between the two cases.
\begin{prop}
\label{noridge}
Let $\tau > 0$ and assume ${\rm{rank}}(\sig) \geq d$. 
Then, $\argmin G_\tau=\emptyset$ if and only if there 
exists $y\in\{1,\dots,h\}$ such that $\sig \Xc \neq 0$. 
\end{prop}
\noindent To solve the optimization problem~(\ref{SIRreg}), Li and Yin~\cite{liyin07} proposed an alternating least-squares algorithm. At iteration $k+1$, given $A^{(k)}$,
$C^{(k+1)}$ and $A^{(k+1)}$ are updated as:
\begin{eqnarray*}
C^{(k+1)}_y &=&  \left({A^{(k)}}^T \sig^2 A^{(k)}\right)^{-1} {A^{(k)}}^T \sig \Xc, \; y=1,\dots,h,\\
\mbox{vec}\left(A^{(k+1)}\right) &=& \left\{
\sum_{y=1}^h f_y \left({C_y^{(k+1)}}^T \otimes \sig\right)^T 
                 \left({C_y^{(k+1)}}^T \otimes \sig\right) + \tau I_{pd}
\right\}^{-1} \\
&\times & \sum_{y=1}^h f_y  \left({C_y^{(k+1)}}^T \otimes \sig\right)^T
\Xc.
\end{eqnarray*}
The authors claimed that such an algorithm converges. As a consequence
of Proposition~\ref{mainprop}, it is easily seen that the limit is
always degenerated.
\begin{coro}
\label{coro}
Let $\tau>0$ and denote by $(A^*,C^*)$ the limit of the sequence
$\left(A^{(k)},C^{(k)}\right)_k$. Necessarily, $A^*$ is the zero $p\times d$
matrix.
\end{coro}
\noindent In view of this result, the good behavior of this algorithm
on simulated and real data reported in \cite{liyin07}, Section~3 and Section~4
cannot be justified from a theoretical point of view.

\section{An alternative ridge SIR estimator}
\label{secnew}

It is possible to modify the criterion $G$ as follows
\begin{equation}
\label{newcrit}
H_{\tau}(A,C) = G(A,C) +  \tau  \sum_{y=1}^h f_y \|A C_y\|^2.
\end{equation}
The first advantage of $H_\tau$ is to be invariant with respect to 
bijective transformations, {\it i.e.}
$$
H_{\tau}(  A M , M^{-1}C ) = H_\tau(A, C), 
$$
for all regular $d\times d$ matrix $M$. This property is natural since
${\mbox{span}}(MA)={\mbox{span}}(A)$. Second, it is readily seen that
the minimization of $H_{\tau}$ does not require the existence of
$\sig^{-1}$ since $H_{\tau}$ can be rewritten as
$$
H_{\tau}(A,C) - H_\tau(0,0) =  \sum_{y=1}^h f_y C_y^T A^T (\sig+\tau I_p) A C_y
- 2  \sum_{y=1}^h f_y \Xc^T A C_y.
$$
Finally, remarking that the original criterion $G$ of SIR~(\ref{SIR})
can also be expanded as
$$
G(A,C) - G(0,0) =  \sum_{y=1}^h f_y C_y^T A^T \sig A C_y
- 2  \sum_{y=1}^h f_y \Xc^T A C_y,
$$
it appears that $H_\tau(A,C)-H_\tau(0,0)$ can be deduced from 
$G(A,C)-G(0,0)$ by substituting $\sig+\tau I_p$ to $\sig$.
Consequently, the estimator of $A$ obtained by minimizing~(\ref{newcrit})
is the Regularized SIR estimator introduced in~\cite{Zhong} since
its columns are the eigenvectors of
$(\sig+\tau I_p)^{-1}\hat\Gamma$ associated to the $d$ largest eigenvalues.
As a conclusion, the introduction of the new functional~(\ref{newcrit})
provides a theoretical framework for the Regularized SIR estimator~\cite{Zhong}.
Thus, a crossvalidation criterion could be derived, similarly to~(8)
in~\cite{liyin07}, for selecting the regularization parameter~$\tau$.

\newpage
\section*{Appendix}

\noindent {\bf{Proof of Proposition \ref{mainprop}}} $-$ 
Let us remark that
\begin{eqnarray}
\nonumber
G_{\tau}(A,C) & = & \sum_{y=1}^h f_y \left( \|\Xcc\|^2 -2 \Xc^T \sig A C_y  + C_y^T A^T \sig^2 A C_y \right ) \\
\label{critere}
&+& \tau \| \vecA \|^2. 
\end{eqnarray}
Using the equality (see for instance~\cite{Harv}, Chapter~16, equation~(2.13)),
\begin{equation} 
\label{trans}
\sig A C_y = \ctsig \vecA,
\end{equation}
for all $y=1,\ldots,h$ and denoting ${\tilde{a}}=\vecA$, we thus have:
\begin{eqnarray*}
G_{\tau}(A,C) = G_{\tau}^*({\tilde{a}},C) & = & \sum_{y=1}^h f_y
 \left\{ \|\Xcc\|^2 
  -  2 \Xc^T \ctsig {\tilde{a}} 
\right.\\
 \ & + & \left . {\tilde{a}}^T \ctsig^T \ctsig {\tilde{a}} \right\} + \tau \| {\tilde{a}} \|^2. 
\end{eqnarray*}
Suppose $\argmin G_\tau\neq \emptyset$ and consider 
$$
({\hat{A}},{\hat{C}})\in\argmin_{A,C} G_{\tau}(A,C).
$$
From~\cite{magneu88}, pp. 119-120, it follows that,
necessarily, $({\hat{A}},{\hat{C}})$ is a stationary point of $G_\tau$
and thus satisfy the set of equations:
\begin{equation}
\label{systgrad}
\grad_i G_{\tau}^*({\hat{a}},{\hat{C}}_1,\ldots,{\hat{C}}_h) = 0, \ i=1,\ldots,h+1, 
\end{equation}
where ${\hat{a}}={\rm{Vec}}({\hat{A}})$, ${\hat{C}} = ({\hat{C}}_1,\ldots,{\hat{C}}_h)$ and 
$\grad_i$ denotes the gradient of $G_{\tau}^*$ 
with respect to its $i$th argument,  $i=1,\ldots,h+1$.
 Straightforward calculations lead to:
\begin{eqnarray}
\nonumber
\grad_1 G_{\tau}^*({\hat{a}},{\hat{C}}_1,\ldots,{\hat{C}}_h) & = & 2 \sum\limits_{y=1}^h  f_y \left\{ \btsig^T \btsig {\hat{a}} \right . \\
\label{eq1}
 \ & - & \left . \btsig^T \Xc \right\} + 2\tau {\hat{a}}, 
\end{eqnarray}
and, for $y=1,\ldots,h$,
\begin{eqnarray}
\nonumber
 \grad_{y+1}G_{\tau}^*({\hat{a}},{\hat{C}}_1,\ldots,{\hat{C}}_h) & = & \grad_{y+1}G_{\tau}({\hat{A}},{\hat{C}}_1,\ldots,{\hat{C}}_h) \\
\label{eq0}
 \ & = & 2 f_y \left\{ {\hat{A}}^T \sig^2 {\hat{A}} {\hat{C}}_y - {\hat{A}}^T \sig \Xc \right\}.
\end{eqnarray}
Thus, multiplying of the left by ${\hat{C}}_y^T$ and using (\ref{trans}),
it follows 
\begin{eqnarray}
\nonumber
 {\hat{C}}_y^T \grad_{y+1}G_{\tau}^*({\hat{a}},{\hat{C}}_1,\ldots,{\hat{C}}_h)
 & = & 2 f_y  \left\{ {\hat{C}}_y^T ({\hat{A}}^T \sig^2 {\hat{A}}) {\hat{C}}_y - {\hat{C}}_y^T {\hat{A}}^T \sig \Xc \right\} \\
\nonumber
 \ & = & 2 f_y \left\{ {\hat{a}}^T \btsig^T \btsig {\hat{a}} \right . \\
\label{eq2}
 \ & - & \left . {\hat{a}}^T \btsig^T \Xc \right\}. 
\end{eqnarray}
Hence, collecting~(\ref{eq1}) and~(\ref{eq2}), it appears that
\[ {\hat{a}}^T \grad_1 G_{\tau}^*({\hat{a}},{\hat{C}}_1,\ldots,{\hat{C}}_h) = \sum_{y=1}^h {\hat{C}}_y^T \grad_{y+1}G_{\tau}^*({\hat{a}},{\hat{C}}_1,\ldots,{\hat{C}}_h) + 2 \tau \| {\hat{a}} \|^2. \]
Since the regularization parameter $\tau$ is positive, condition~(\ref{systgrad})
implies $\| {\hat{a}} \|^2=0$, {\it i.e.}  ${\hat{A}}$ is the zero $p\times d$
matrix. Replacing in~(\ref{critere}), we have
\[ G_{\tau}({\hat{A}},C) = G_{\tau}(0,C) = \sum_{y=1}^h f_y \|\Xcc\|^2, \]
for all $C \in \R^{d \times h}$ and the result is proved.
\CQFD 
 \ \\
\noindent {\bf{Proof of Corollary \ref{coro}}} $-$ The limit of the
sequence verifies the set of equations
\begin{eqnarray*}
C^*_y &=&  ({A^*}^T \sig^2 A^*)^{-1} {A^*}^T \sig \Xc, \; y=1,\dots,h,\\
\mbox{vec}\left(A^*\right) &=& 
 \left\{ \sum_{y=1}^h f_y \left({C_y^{*}}^T \otimes \sig\right)^T 
                 \left({C_y^{*}}^T \otimes \sig\right) + \tau I_{pd}
\right\}^{-1} \\
&\times & \sum_{y=1}^h f_y  \left({C_y^{*}}^T \otimes \sig\right)^T
\Xc.
\end{eqnarray*}
Thus, from~(\ref{eq0}) it follows that
$$
\grad_{y+1}G_{\tau}(A^*,C^*_1,\ldots,C^*_h)=0
$$
for all $y=1,\dots,h$, while, from~(\ref{eq1}),
$$
\grad_{1}G_{\tau}(A^*,C^*_1,\ldots,C^*_h)=0.
$$
Consequently, $(A^*,C^*)$ is a stationary point of $G_\tau$, and,
following the proof of Proposition~\ref{mainprop}, necessarily $A^*$
is the zero $p\times d$ matrix.\CQFD

\noindent {\bf{Proof of Proposition \ref{noridge}}} $-$ 
First, let us suppose that $\sig \Xc =0$ for all $y\in\{1,\dots,h\}$. Then,
\begin{eqnarray*} 
G_{\tau}(A,C) & = & \sum_{y=1}^h f_y \|\Xcc\|^2 + \sum_{y=1}^h f_y C_y^T (A^T \sig^2 A) C_y  + \tau \| \vecA \|^2\\
&\geq & \sum_{y=1}^h f_y \|\Xcc\|^2 \\
&=& G_{\tau}(0,C),
\end{eqnarray*}
which entails that $G_{\tau}(A,C)$ is minimum for every $C$ if $A$ is the zero matrix. As a consequence $\argmin G_\tau\neq \emptyset$.
This concludes the first part of the proof. 
Conversely, suppose there exists $y_0 \in \{ 1, \ldots,h\}$ such that $\sig \Xy \neq 0$. 
Let $\tau>0$ and let us 
prove that there exist $A \in \R^{p \times d}$ and $C \in \R^{d \times h}$ such that $G_{\tau}(A,C) < G_\tau(0,C)$. 
To this end, let $q_i$, $i=1,\ldots,p$ be the eigenvectors of $\sig$ associated to  the eigenvalues $\lambda_i$, $i=1,\ldots,p$. Since
\[ \sig \Xy = \sum_{i=1}^p \lambda_i q_i q_i^T \Xy \neq 0, \]
there exists an eigenvector $q^*$ associated to a random value $\lambda^*>0$ such that $q^*{q^*}^T \Xy \neq 0$. Thus $\|\Xy^T q^* \| \neq 0$,  
and let $\varepsilon$ such that:
\begin{equation}
\label{eps}
0 < \varepsilon < \sqrt{\frac{f_{y_0}}{\tau d} }\| \Xy^T q^* \|. 
\end{equation}
The matrices $A$ and $C$ are defined as follows. 
The first column of $A$ is the vector $\varepsilon q^*$ 
and the $d-1$ following columns of $A$ are the vectors $\varepsilon q_{j_i}$, $i=1,\ldots,d-1$ where the $q_{j_i}$'s are orthogonal eigenvectors (with unit norm) of $\sig$ associated to positive eigenvalues $\lambda_{j_i}$'s. 
Note that, since ${\rm{rank}}(\sig) \geq d$, such a matrix $A$ always exists. 
All the columns of $C$ are chosen to be the null vector except the $y_0$th
one defined by:
\[ C_{y_0} = (A^T \sig^2 A)^{-1} A^T \sig \Xy = \frac{1}{\varepsilon} \left ( \frac{q^*}{\lambda^*},\frac{q_{j_1}}{\lambda_{j_1}},\ldots, \frac{q_{j_{d-1}}}{\lambda_{j_{d-1}}} \right )^T \Xy. \]
Such choices entail
\begin{eqnarray*}
&&G_\tau(A,C)-G_\tau(0,C)\\
&=&\sum\limits_{y=1}^h f_y \left\{ C_y^T (A^T \sig^2 A) C_y - 2 \Xc^T \sig A C_y \right \} + \tau \| \vecA \|^2 \\
& = & f_{y_0} \left\{ C_{y_0}^T (A^T \sig^2 A) C_{y_0} - 2 \Xy^T \sig A C_{y_0} \right\} + \tau \| \vecA \|^2 \\
& = & -f_{y_0} \Xy^T \sig A (A^T \sig^2 A)^{-1} A^T \sig \Xy + \tau \| \vecA \|^2 \\
& = & -f_{y_0} \| \Xy^T q^* \|^2 -f_{y_0} \sum\limits_{i=1}^{d-1} \| \Xy^T q_{j_i} \|^2 + \tau d \varepsilon^2 \\
& \leq & -f_{y_0} \| \Xy^T q^* \|^2 + \tau d \varepsilon^2\\
& < & 0,
\end{eqnarray*}
from (\ref{eps}). Thus, $(0,C)\notin \argmin G_\tau$ and taking account of
Proposition~\ref{mainprop} yields $\argmin G_\tau= \emptyset$. 
\CQFD

\end{document}